\documentclass[twocolumn,10pt,journal]{asme2e}

\usepackage{graphicx} 
\usepackage{amsmath}
\usepackage{cancel}
\usepackage{graphicx}

\title{Toward Fast Topological-Shape Optimization With Boundary Elements}

\author{Igor Ostanin\thanks{Corresponding author, tel: +79150174677, e-mail:i.ostanin@skoltech.ru}
    \affiliation{
	Skolkovo Institute \\
      of Science and Technology\\
	100 Novaya St., Skolkovo\\
      143025, Russia\\
    }
}

\author{Denis Zorin
    \affiliation{ 
	Courant Institute of Mathematical Sciences \\
     New York University, \\
     251 Mercer St, New York\\
     NY 10012, United States\\
     Skolkovo Institute \\
     of Science and Technology\\
    100 Novaya St., Skolkovo\\
    143025, Russia\\
    }
}

\author{Ivan Oseledets
        \affiliation{
	Skolkovo Institute \\
      of Science and Technology\\
	100 Novaya St., Skolkovo\\
      143025, Russia\\
    }
}

\IfFileExists{revision.tex} {
    \input{revision}
} {
    
}

\begin{document}

\maketitle

\begin{abstract}
 {\it Abstract. Wide variety of engineering design tasks can be formulated as constrained optimization problems where the shape and topology of the domain are optimized to reduce costs while satisfying certain constraints. Several mathematical approaches were developed to address the problem of finding optimal design of an engineered structure. Recent works \cite{BEM2DMarczak, BEM3D} have demonstrated the feasibility of boundary element method as a tool for topological-shape optimization. However, it was noted that the approach has certain drawbacks, and in particular high computational cost of the iterative optimization process. In this short note we suggest ways to address critical limitations of boundary element method as a tool for topological-shape optimization. We validate our approaches by supplementing the existing complex variables boundary element code for elastostatic problems with robust tools for fast topological-shape optimization. The efficiency of the approach is illustrated with a numerical example.
 }
\end{abstract}

\section{Introduction}

Problems of optimal design, \textit{i.e.} variation of shape and topology
of the domain to extremize certain functional subject to additional
constraints, are ubiquitous in different branches of engineering.
Recent emergence and rapid development of stereolithography and 3D
printing technologies \cite{3Dprint,Stereo} have enabled fast prototyping
of complex-shaped structures and revived the interest in automated
optimal design. The most common formulation of optimization problems
studied in structural engineering (e.g. \cite{Novotny3Delast}) seeks
for an optimal shape and topology of an elastic body that minimize
the strain energy (compliance) while satisfying the weight constraint
and the additional constraints imposed by the boundary value problem.
Non-convex nature of such optimization problems often makes finding
globally optimal designs difficult or impossible. Various numerical
methods, including shape gradient-based approaches \cite{gradient_based},
level set methods \cite{AllaireLevelSet,LevelSet}, homogenization
\cite{Homog2DKikuchi,AllaireHomog} and topological-shape sensitivity
\cite{NovotnySokolovskyBook} have been developed during the last
few decades to address this class of problems. These approaches are
typically implemented within the finite element method (FEM) context.
In case if admissible designs include only homogeneous regions with
piecewise-constant properties (micro-structured composite designs
are prohibited), the problem of optimization reduces to finding optimal
configuration of the domain boundaries. Few recent papers \cite{BEM2DMarczak,BEM3D}
suggested that boundary integral approaches, and in particular boundary
element method (BEM), can be a convenient tool to address this class
of problems. The implementation of optimization algorithms within
BEM utilizes the concept of topological derivative (TD) \cite{Sokolmain,NovotnySokolovskyBook,Sokol3D}
- a cost of making an infinitesimal circular (spherical) cavity with
a center in a given point of the domain. Existence of this kind of
derivative opens wide avenue for numerous gradient-based approaches.
The most straightforward one utilizes the idea of calculation of the
field of TD on a given mesh, and removing material in the regions
where the TD is below certain cutoff level. This kind of approach
has demonstrated its feasibility in both 2D and 3D elastostatic problems
\cite{BEM2DMarczak,BEM3D}. However, existing papers on the subject
do not address the issue of complexity of topological-shape optimization
with BEM, and, particularly, the question of competitiveness of BEM
approaches as compared to FEM. 

In this note we suggest a simple optimization technique based on 2D
BEM. It is demonstrated that even without using fast BEM techniques,
one can easily reach $O(n_{s}^{2})$ asymptotic performance of shape
optimization with BEM, which corresponds to the performance of 2D
FEM approaches - $O(n_{v})$ ($n_{s}$ and $n_{v}$ are the numbers
of surface elements and volume elements after BEM/FEM discretization). 

The suggested technique is based on complex variables boundary element
method (CVBEM) \cite{CVBEM_main}. It utilizes quadtree-based calculation
of TD inside the domain, allowing $O(n_{s}^{2})$ performance of evaluation
of fields inside the domain. We also use simple $O(n_{s}^{2})$ algebraic
solver based on blockwise update of the inverse system matrix to solve
direct boundary value problem at every iteration. 

The paper has the following organization. The next section describes
the method of topological optimization employed in our research, as
well as its novel features. Section 3 provides an illustrative numerical
example that is served to highlight the capabilities of our approach.
The summary and discussion of our results, as well as the major directions
of future work are presented in Section 4.

\section{Method}

\subsection{Optimization technique}

In this paper we develop our approach in 2D using the framework of
CVBEM \cite{CVBEM_main}. The method offers a number of attractive
features: i) complex hypersingular boundary integral equation (BIE)
for tractions and displacement discontinuities, derived for a system
that can include an infinite matrix, finite-sized blocks with piecewise-constant
properties, voids and cracks \cite{CVBEM_main}; ii) closed form analytical
calculation of all the integrals in BIE \cite{CVBEM_main}, iii) circular
elements to approximate curved boundaries \cite{CVBEM_main}; vi)
asymptotic approximations to model cracks \cite{CVBEM_main}; v) implementation
of piecewise-constant body forces \cite{cvbem_bf}. 

\begin{figure}
\includegraphics[width=9cm]{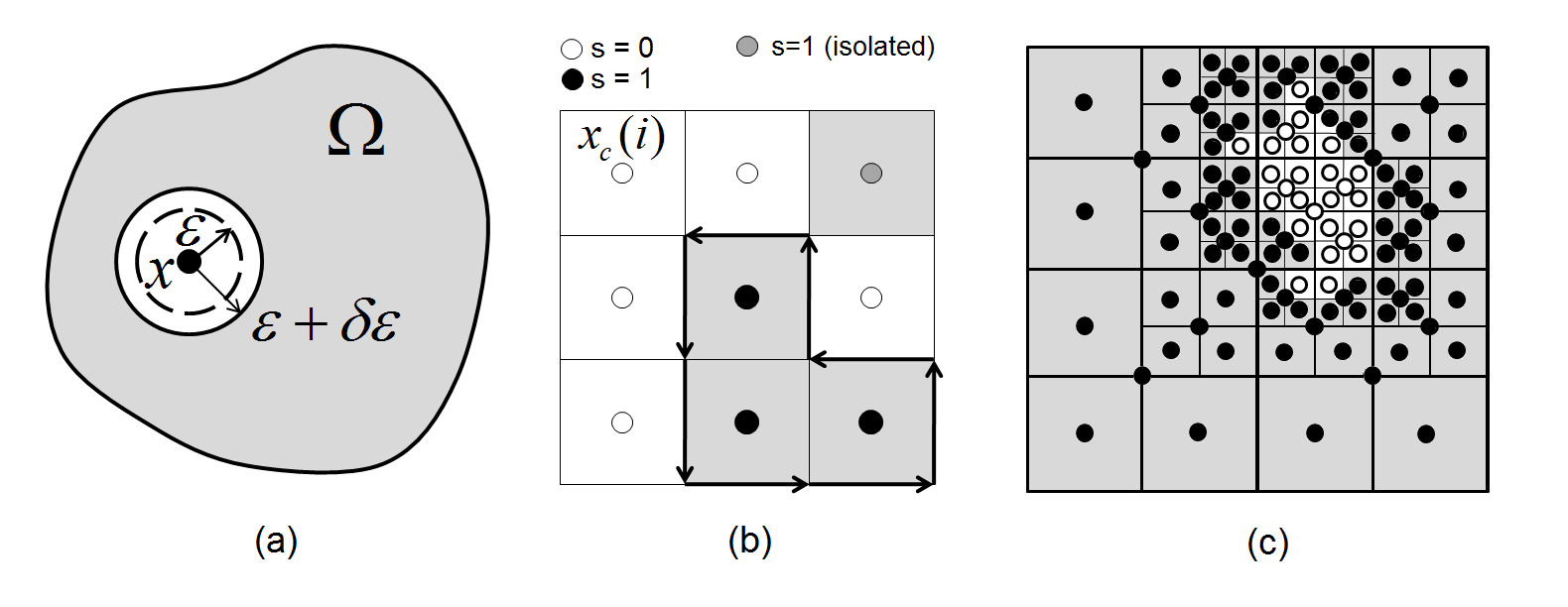}

\protect\caption{(a) Details on the definition of TD. (b) Optimization on the rectangular
mesh. (c) Quadtree strategy of sampling of the TD inside the domain. }
\end{figure}

Our scheme of topological optimization utilizes the notion of topological
derivative \cite{NovotnySokolovskyBook}. Following \cite{Novotny2Delast}
we define the TD as (Fig. 1(a)):

\begin{equation}
D(x)=\underset{\underset{\delta\varepsilon\rightarrow0}{\varepsilon\rightarrow0}}{\lim}\frac{\varPsi(\varOmega_{\varepsilon+\delta\varepsilon})-\varPsi(\varOmega_{\varepsilon})}{f(\varepsilon+\delta\varepsilon)-f(\varepsilon)}
\end{equation}

Here $\Psi$ is the cost functional, $\Omega_{\varepsilon}$ is the
original domain perturbed by the presence of an infinitesimal cavity
of radius $\varepsilon$ centered at the point $x$, $\delta\varepsilon$ is
the small perturbation of the cavity radius, and $f$ is regularizing,
problem-dependent function. Here and below we'll be working with the
strain energy functional, that serves as a global measure of compliance
of an elastic structure. For this case the area (volume) of the cavity
serves as regularizing function. In absence of body forces, the closed
form analytical expressions are available both in 2D and 3D \cite{Novotny2Delast,Novotny3Delast}.
In the case of 2D plane strain elasticity the expression for TD writes

\begin{equation}
D(x)=\frac{2}{(1+\nu)(1-2\nu)}\sigma\cdot\epsilon+\frac{(1-\nu)(4\nu-1)}{2(1-2\nu)}tr\sigma tr\epsilon
\end{equation}

Here $\sigma$ and $\epsilon$ are stress and strain tensors at point
$x$. In case of presence of a constant body force $\boldsymbol{b}$
the expression (2) should be enriched with an additional term, associated
with work done by the body force $-\boldsymbol{bu}(x)$, where $\boldsymbol{u}(x)$
is displacement at point $x$. The treatment of constant body force
within CVBEM is discussed in \cite{cvbem_bf}. It worth noting here
that similar analytical expressions for TDs are also available for
other important classes of cost functionals, in particular, for the
components of homogenized elasticity tensor of a periodic cell \cite{periodicCell1,periodicCell2}.

Shape optimization approaches based on TDs use either hard-kill methods
\cite{Novotny2Delast,BEM2DMarczak}, or bidirectional optimization\cite{BESO}.
In our work we utilize a hard-kill algorithm, in which the optimization
is achieved by progressive elimination of material in the areas where
the TD is below the cutoff level, as discussed further.

\subsection{Boundary generation procedure}

The new set of domain boundaries is created at every iteration of optimization procedure. 
Boundary generation strategy employed in our work is similar to one
described in \cite{BEM3D}. We discretize the initial domain onto
a set of $M$ square cells (Fig. 1(b)). The fields inside the domain
are calculated at the center of each cell $x_{c}(i)$, $i=1..M$.
The cutoff level for the TD is defined as $D_{0}=D_{min}+\alpha(D_{max}-D_{min})$,
where $D_{min}$ and $D_{max}$ are minimum and maximum values of
TD within the current domain, and the coefficient $\alpha$ is tuned
in range between 0.1\% and 2\%. For every cell of the initial domain
we define the Boolean status $s$. At the beginning of the iterative
process $s=1$ for every cell. At every iteration we assign $s=0$
for the cells with $D(x_{c}(i))<D_{0}$ (marked with empty points
in Fig. 1(b)), and for isolated cells (marked with gray points in
Fig. 1(b)). When the status is assigned, we generate the new boundary
using straightforward algorithm - if $i$-th cell has $s=1$ and its
right (top, left, bottom) neighbor has $s=0$, then generate a corresponding
boundary element. For every boundary between the neighboring cells
one straight element with three points of collocation (quadratic approximation)
is used. At every iteration we mark the elements that were deleted
and those that were added at the current step. 

Boundary value problem constraints are incorporated explicitly - for
the cells bounded by the elements with non-homogeneous Neumann boundary
conditions $s\equiv1$. The volume constraint does not present explicity
in this scheme, so it is incorporated as the stopping criteria - the
iterative process is discontinued when the ratio of current area of
the material to the initial one $R$ reaches the prescribed value
$R_{0}$.

\subsection{Quadtree sampling of fields inside the domain}

The described procedure of sampling of TDs on the uniform mesh should
be considered computationally inefficient, since the calculation of
a field of TDs inside the domain takes $O(n_{s}^{3})$ operations.
In order to reduce the asymptotic complexity of this step, we use
the specific strategy of calculation of TDs, which is based on quadtree
algorithm of sampling \cite{qtree}. The main idea is to use few different
levels of refinement, that are employed for initial detection and
further refinement of features of optimized domain. Within this approach
each refined cell is subdivided onto 4 sub-cells (Fig. 1(c)). One
can formulate different possible criteria of refinement. In this work we use
the following criterion - if the values of $s$ for a current cell
and its nearest neighbor are not the same, both cells should be refined
to the next level. The boundary generation algorithm described above
is used to generate boundaries around the points of highest level
of refinement. It is easy to see that such an algorithm of sampling and
boundary generation requires calculation of the topological derivative
at $O(n_{s})$ points inside the domain, leading to $O(n_{s}^{2})$
operations of evaluations of boundary integrals for each boundary
element.

\subsection{Iterative update of the inverse matrix }

CVBEM \cite{CVBEM_main} generates non-symmetric and non-sparce system
matrix, and the corresponding system of equations is difficult to
solve using regular iterative approaches. Here we describe one possible
way to treat the resulting system matrix during the iterative updates
of the initial boundary value problem. The original boundary
consisting of $n_{s}$ elements leads to $N_{s}=6n_{s}$ rows/columns
in the resulting non-symmetric matrix generated by CVBEM (considering
$3$ collocation points per element and $2$ degrees of freedom per
collocation point). Assume that at $k-th$ iteration $n_{a}$
elements have been added to a boundary, and $n_{r}$ elements have
been removed. This results in $N_{a}=6n_{a}$ and $N_{r}=6n_{r}$
added and removed rows/columns in the system matrix. If, as always
the case in the iterative process, $n_{a}\ll n_{s}$, $n_{r}\ll n_{s}$, it
is efficient to perform update of the inverse matrix using incremental
expressions, such as Sherman\textendash Morrison\textendash Woodbury
formula \cite{swm} or Banachiewicz \cite{Banachiewicz} formula for
blockwise matrix inverse:

\begin{equation}
\left(\begin{array}{cc}
A & B\\
C & D
\end{array}\right)^{-1}=\left(\begin{array}{cc}
A^{-1}+A^{-1}BS{}^{-1}CA^{-1} & -A^{-1}BS{}^{-1}\\
-S^{-1}CA^{-1} & S^{-1}
\end{array}\right)
\end{equation}
Where $S=D-CA^{-1}B$. Denoting blocks of an extended matrix as $E,F,G$
and $H$, and expressing $A^{-1}$, we obtain the formula for inverse
matrix update in case of removing elements:

\begin{equation}
A^{-1}=\left(\begin{array}{cc}
E & \cancel{F}\\
\cancel{G} & \cancel{H}
\end{array}\right)^{-1}=E-FH^{-1}G
\end{equation}

It is clear that calculations of both expressions do not require $(N_{s}\times N_{s})(N_{s}\times N_{s})$
matrix multiplications, only lower-rank operations with $(N_{s}\times N_{s})$,
$(N_{s}\times N_{a})$,$(N_{s}\times N_{d})$,$(N_{a}\times N_{a})$
and $(N_{d}\times N_{d})$ matrices. This leads to $O(n_{s}^{2})$
performance of the algorithm of iterative matrix update.

It is important to mention that in case if the matrix update is associated
with the new cavity in the domain, matrix $S$ is singular and requires
regularization, \textit{e.g.} truncated SVD regularization \cite{SVD_regularization}
(used in this work) or explicit incorporation of additional conditions
that fix rigid body motion of a new closed boundary \cite{CVBEM_main}.

The suggested (or similar) schemes of low-rank inverse matrix update
are well know and are widely used in adjacent fields of scientific
and engineering calculations, including network structures, asymptotic
analysis and boundary value problems with changing boundaries (see
the excellent review presented in \cite{Inverse}). It is thus natural
to adopt this approach for fast iterative optimization schemes with
boundary elements.

\section{Numerical example}

In this section we discuss a simple benchmark example that allows
us to evaluate the capabilities of our approach. Consider a 2D problem
of optimization of a shape and topology of a fixed support. The initial
domain has the square shape (Fig. 2(a)). The left side of the support
is rigidly fixed, and the point load is applied at the upper-right
angle.

\begin{figure}
\includegraphics[width=9cm]{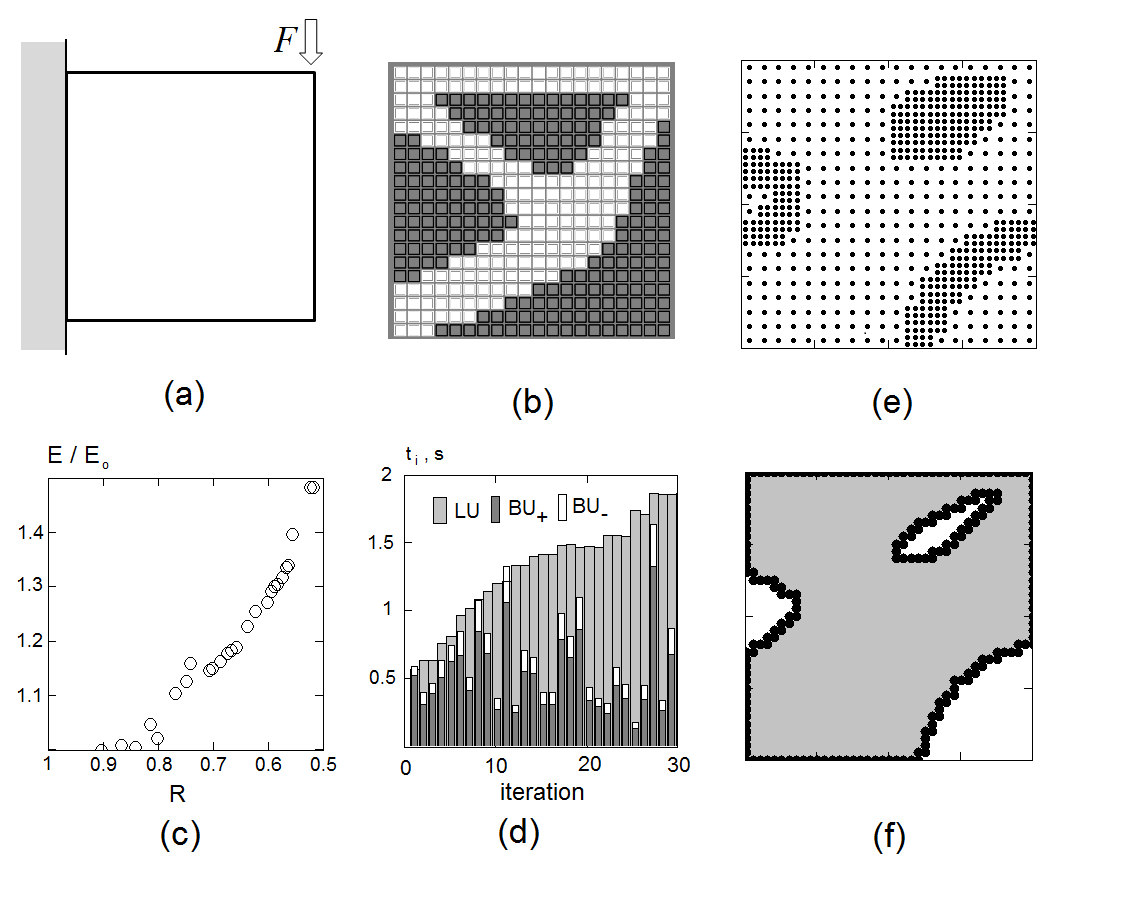}

\protect\caption{Shape optimization of a fixed support. (a) Initial boundary value
problem. (b) Final configuration reached after $29$ iterations. (c)
The ratio of current strain energy $E$ to the initial strain energy
$E_{0}$ as a function of current volume fraction $R$.}
\end{figure}

We utilize the optimization scheme described above to find an optimal
shape and topology of the support, providing the smallest strain energy
$E$ for the volume fraction $R_{0}=0.5$. The cutoff parameter $\alpha$
was set to $0.003$. The TDs are sampled on a grid containing $20\times20$
cells. The solution was found in $29$ iterations. The final shape
is given in Fig. 2(b). The obtained solution is in good agreement
with the one obtained with FEM and BEM in the earlier works \cite{Novotny2Delast,BEM2DMarczak}.
Fig. 2(c) gives the evolution of the strain energy functional during
the iterative process. Fig. 2(d) presents the comparison of times
spent during the iteration with full LU factorization ($LU$), and
the blockwise inverse matrix update, including adding rows/columns
according to (3) ($BU_{+}$) and removing rows/columns according to
(4) ($BU_{-}$). The boost in performance clearly depends on rank
of the update (number of added/removed elements) and varies significantly
from iteration to iteration. The total time spent on the full LU factorization
is 41 s, whereas the total time spent on the blockwise matrix update
is 15.5 s. Calculations were performed with regular Core i-5 laptop
machine. Linear algebra operations were implemented within non-optimized
BLAS/LAPACK framework \cite{lapack} and win-32 gfortran compiler
\cite{Fortran}.

In order to illustrate the quadtree sampling of TDs inside the domain,
we consider the same example, but with the twice higher degree of
grid refinement. Two levels of grid refinement are used. The coarse
level discretizes domain onto $20\times20$ cells. As we could see
above, this level is sufficient to detect the important features of
the optimized domain. On the finer level we discretize the domain
onto $40\times40$ cells. This level is used to render finer features
of the optimal design. Fig. 2(e) gives the sampling points during
first iteration, Fig. 2(f) demonstrates corresponding configuration
of domain boundaries. Using quadtree sampling has led to decrease
of the number of points inside the domain from 1600 to 676, and decrease
of the corresponding computational time from 11.3 to 4.8 s.

\section{Discussion and conclusions}

In this work we suggested a set of tools for topological-shape optimization
with boundary elements. As we could see on the example of 2D CVBEM
method and a simple toolkit for topological optimization, one can
reach the computational performance available with 2D FEM techniques.

It is therefore clear that the usage of fast BEM techniques \cite{fastBEM,fastBEM2}
would undoubtedly allow to outperform all the existing FEM techniques
of topological-shape optimization. For example, using fast BEM would
allow solution of direct boundary value problem for $O(n_{s})$ operations,
and calculation of the field of TD inside the domain for $O(n_{s}ln(n_{s}))$
operations, which is much faster than the corresponding operations
performed within FEM (both take $O(n_{v})$ operations). In addition,
FEM techniques require good quality domain discretization, generation
of which takes at least $O(n_{v}ln(n_{v}))$ operations.

These considerations motivate the development of three-dimensional
fast BEM-based framework for topological-shape optimization of an
elastic domain. Its implementation can be based on the principles
similar to those presented in this work, and extended with additional
features: smoother boundaries generation, combined shape and topology
optimization iterations, \textit{etc}. Clearly, the BEM-based techniques
should become the most computationally efficient tools for topological-shape
optimization.

\begin{acknowledgment}
\end{acknowledgment}

%
% ...

\bibliographystyle{asmems4}

\bibliography{ostanin}


\end{document}